\documentclass[12pt]{article}
\usepackage{amsmath, amssymb}

\begin{document}

\title{
Further Comments on Realization of Riemann Hypothesis via Coupling Constant Spectrum}

\author{R. Acharya \\
Department of Physics\\
Arizona State University\\
Tempe, AZ 85287-1504}

\maketitle

\begin{abstract}
We invoke \textbf{Carlson's theorem} to justify \textbf{and} to
confirm the results previously obtained on the validity of Riemann
Hypothesis via the coupling constant spectrum of the zero energy
$S$-wave Jost function a la N.~N.~Khuri, for the real, repulsive
inverse-square potential in non-relativistic quantum mechanics in
3 dimensions.
\end{abstract}

In a previous note (hereinafter referred to as I) \cite{Acharya},
we have exhibited the solution for the $S$-wave Jost function at
zero energy (i.e., $k^{2}=0$) for the real, repulsive
inverse-square potential, $V(r)=\frac{\lambda}{r^{2}}$ ($\lambda >
0$). We then argued that it is justified to identify the $S$-wave
Jost function, $\mathcal{X}(s)\equiv F_{+}((s(s-1)$; $k^{2}=0)$ with
the Riemann $\xi(s)$ function \cite{Borwein} (up to an
entire function of order one, with \textbf{No} zeros) where the
coupling constant $\lambda$ is \textbf{constrained} to satisfy the
relation, i.e.,
\begin{equation}\label{eq:1}
\lambda = s(s-1)
\end{equation}

The argument hinged on demonstrating that the $S$-wave Jost
function has \textbf{only} zeros on the line $s_{n} = \frac{1}{2}
+ i \gamma^{'}_{n}$, $n = 1$,~2,~3, \dots $\infty$. This was
accomplished by first demonstrating that the analogue of Eq.~(1.2)
of Khuri's paper \cite{Khuri} (\textbf{without} the presence of the potential
$V(r)$ under the sign of integration), i.e.,
\begin{equation}\label{eq:2}
\big[
    Im\lambda_{n}(i\tau)
\big]
    \int^{\infty}_{0}
\big|
    f(\lambda_{n}(i\tau); \text{ } i\tau; \text{ } r
\big| ^{2} dr = 0
\end{equation}
holds, where the all-important constraint, i.e.,
\begin{equation}\label{eq:3}
    \int^{\infty}_{0}
\big|
    f(\lambda_{n}(i\tau); \text{ } i\tau; \text{ } r
\big|
    ^{2} dr < \infty
\end{equation}
is satisfied, i.e.,
\begin{equation}\label{eq:4}
\big[
    Im\lambda_{n}(i \tau)
\big]
    \frac{1}{\tau} \frac{2}{\pi} \int^{\infty}_{0} r K^{2}_{\nu(i
    \tau)} (r) dr = 0
\end{equation}
where
\begin{equation}\label{eq:5}
\int^{\infty}_{0} r K^{2}_{\nu (i \tau)} (r) dr = \frac{1}{8}
\frac{\pi \nu (i \tau)}{\sin \pi \nu (i \tau)}
\end{equation}
with $\nu (i \tau) = \sqrt{\lambda (i \tau) + \frac{1}{4}}$,
$\nu (i \tau) \neq 1$,~2,~3, \dots $\infty$.
Then, Eq.~(\ref{eq:3}) follows for $\tau \neq 0$ ($\tau > 0$).

Thus,
\begin{equation}
I m \lambda_{n} (i\tau) = 0
\end{equation}
We then followed Khuri and arrived at the conclusion,
\begin{equation}
\lambda_{n} (0) \equiv s_{n} (s_{n} -1)
\end{equation}
for \textbf{all} {\mathversion{bold}$n$}, is \textbf{real} and
\textbf{negative}. When $\lambda_{n}(0)$, for all $n$ is real and
negative, (and by rescaling such that $\lambda_{n} <
-\frac{1}{4}$, $\nu(0) = \sqrt{\lambda_{n}(0)+\frac{1}{4}}$
becomes imaginary and hence R.H.S. of Eq.~(\ref{eq:5}) will become
the ratio of $\frac{\pi \nu(0)}{\sinh \pi \nu(0)}$, at $\tau = 0$.
In other words, the zero energy coupling spectrum lies on the
negative real line for $V(r) = \lambda V^{*}(r)$, $V(r) =
\frac{1}{r^{2}}$, $V^{*} \geqslant 0$. In view of
Eq.~(\ref{eq:1}), i.e., $\lambda = s(s-1)$, one then concluded
that $\mathcal{X}(s)$ has all its infinite number of nonzero zeros
on the ``critical" line, Re $s_{n} = \frac {1}{2}$.

The \textbf{key question} arises whether or not, the positions
$\gamma_{n}$ of the infinite number of non-zero zeros on the
critical line, \textbf{coincide} with the Riemann Hypothesis,
i.e.,

$$
s_{n} = \frac{1}{2} \pm i \gamma_{n}, \text{ } n=1 \text{,~2,~3
\dots} \infty
$$
(Recall that Hardy has proved that $\xi(s)$ has \textbf{infinite number}
of zeros on the critical line.)

In other words, one \textbf{must} demonstrate that the positions
of the zeros obtained from the potential, $V(r) =
\frac{\lambda}{r^{2}}$, i.e.
$$
s^{'}_{n} = \frac{1}{2} \pm i \gamma^{'}_{n}, \text{ } n=1\text{,~2,~3
\dots} \infty
$$
coincide exactly with the positions ($\gamma_{n}$) of the Riemann
zeros. This was \textbf{not} done in our earlier note! \textbf{We propose
to demonstrate this here.} Our tool is Carlson's theorem
\cite{Levin}.

Recall,
\begin{equation}
\mathcal{X}(s^{'}) = e^{\pi s^{'}}.
\end{equation}
and
\begin{equation}
\xi(s^{'}) = C e^{A(s^{'}+1)}
\end{equation}
for real $s^{'} \geqslant 0$, $s^{'} = s-1$, where we are stating
the \textbf{established} fact that \textbf{both} $\mathcal{X}(s)$
and $\xi(s)$ are \textbf{entire} with \textbf{No zeros} for real
$s^{'} > 0$.

Then, we invoke Hadamard's factorization theorem to set
\begin{equation}\label{eq:10}
e^{a+\frac{bs^{'}}{m}} \mathcal{X} (\frac{s^{'}}{m}) - \xi (\frac{s^{'}}{m}) =
0, \quad \frac{s^{'}}{m} = 1, 2, 3 \dots
\end{equation}
Eq.~({\ref{eq:10}) is satisfied iff
\begin{equation}\label{eq:11}
C e^{A} = e^{a} \text{ and } A = b + \pi
\end{equation}

We now \textbf{invoke Carlson's theorem} \cite{Levin}:\\
In order that each entire function $f(z)$ satisfying the
conditions
\begin{align}
f(z) & = 0(1) e^{\alpha|z|} \text{, for some } \alpha < \infty \\
f(iy) & = 0(1) e^{\beta|y|} \text{, for some } \beta < \pi \\
f(n) & = 0 \text{ for each positive integer `$n$'} \\
\text{then } f(z) & \equiv 0 \text{ !}
\end{align}
Thus, we conclude that
\begin{equation}\label{eq:16}
e^{a+\frac{bs^{'}}{m}} \mathcal{X} \bigg( \frac{s^{'}}{m} \bigg) = \xi \bigg( \frac{s^{'}}{m}
\bigg) \text{, \textbf{all} } s^{'}
\end{equation}

Clearly Eq.~(\ref{eq:10}) is satisfied for \textbf{all} $\frac{s^{'}}{m} >
0$, when Eq.~(\ref{eq:11}) holds. An analytic (entire) function
which vanishes in a subdomain of its region of analyticity, \textbf{must}
necessarily vanish in its ``entire" domain of analyticity, i.e.,
for \textbf{all} $s^{'}$! Thus, in a sense, one need not invoke
Carlson's theorem to arrive at Eq.~(\ref{eq:16})!

It is now obvious that the zeros on the critical line of the zero
energy Jost function $\mathcal{X}(s)$ and the zeros (infinite
number, established by G.~H.~Hardy) of Riemann's $\xi(s)$ function
must necessarily coincide, i.e.,
$$
s^{'}_{n} = \frac{1}{2} \pm i \gamma^{'}_{n} = s_{n} = \frac{1}{2}
\pm i \gamma_{n}, \quad n=1,2, \dots \infty.
$$

We can explicitly verify this assertion as follows:

We invoke the Hadamard
factorization for every entire function $g(z)$ of order one and
``infinite type" (which guarantees the existance of \textbf{infinitely
many} non-zero zeros [valid for \textbf{all} $z$!] \cite{Levin}:
\begin{equation}\label{eq:17}
g(z) = z^{m} e^{B} e^{Dz} \prod^{\infty}_{n=1} \bigg( 1 -
\frac{z}{z_{n}} \bigg) \exp \bigg( \frac{z}{z_{n}} \bigg)
\end{equation}
where `$m$' is the multiplicity of the zeros (so that $m = 0$ for
a simple zero). We can now apply Eq.~(\ref{eq:17}) to the zero
energy Jost function $\mathcal{X}(s^{'})$ and the Riemann's
$\xi(s^{'})$ function, assuming, $s^{'}_{n} = \frac{1}{2} + i
\gamma^{'}_{n}$ for $\mathcal{X}(s^{'})$ and $s^{'}_{n} =
\frac{1}{2} + i \gamma_{n}$ for Riemann's $\xi(s^{'})$.

This leads to the equality [Eq.~(\ref{eq:16})]:
\begin{equation}\label{eq:18}
\begin{aligned}
e^{H} e^{L\frac{s^{'}}{m}} & \prod^{\infty}_{n=1}
\bigg( 1 - \frac{s^{'}}{ms^{'}_{n}} \bigg)
\exp \bigg( \frac{s^{'}}{ms^{'}_{n}} \bigg) \\
 & = e^{J} e^{K\frac{s^{'}}{m}} \prod^{\infty}_{n=1}
\bigg( 1 - \frac{s^{'}}{ms_{n}} \bigg)
\exp \bigg( \frac{s^{'}}{ms_{n}} \bigg).
\end{aligned}
\end{equation}

First set $s^{'} = 0$ in Eq.~(\ref{eq:18}):
\begin{equation}
\Longrightarrow \qquad \qquad H = J
\end{equation}

Next set, {\mathversion{bold}$s^{'} = m s_{n}$}, $s_{n} = \frac{1}{2} + i
\gamma_{n}$
\begin{equation}\label{eq:20}
\begin{aligned}
\Longrightarrow \qquad e^{Ls_{n}} & \prod^{\infty}_{n=1}
\bigg( 1 - \frac{s_{n}}{s^{'}_{n}} \bigg)
\exp \bigg( \frac{s_{n}}{s^{'}_{n}} \bigg) \\
&= e^{Ks_{n}} \prod^{\infty}_{n=1} (1 - 1) \exp(1).
\end{aligned}
\end{equation}

Since RHS of Eq.~(\ref{eq:20}) \textbf{vanishes}, we conclude that
$s_{n}$ must equal $s^{'}_{n}$:
\begin{equation}
\therefore \gamma_{n} = \gamma^{'}_{n}
\end{equation}

This \textbf{establishes the key result} that we set out to prove:
\textbf{The Riemann zeros on the critical line do coincide precisely with
the zeros of the zero energy} {\mathversion{bold}$S$} \textbf{wave Jost function for}
{\mathversion{bold}$V(r) = \frac{\lambda}{r^{2}} \text{ } (\lambda >
0)$} \cite{Footnote}.

\renewcommand\refname{References and Footnotes:}

\begin{center}
\noindent \small{\textbf{Acknowledgement}}
\end{center}

\vspace{-12pt}
I am grateful to Irina Long for
her generous help in the preparation of this manuscript.


\begin{thebibliography}{99}
\bibitem{Footnote} We hasten to correct misstatements in Reference
1. The sentence immediately following Eq.~(30) must be deleted.
Eq.~(40) and Eq.~(41) must be replaced by
\begin{equation}
\alpha = \ln 2, \tag{40}
\end{equation}
\begin{equation}
\beta = 0 \tag{41}
\end{equation}
One arrives at this by setting $\mathcal{X}(0) = \mathcal{X}(1) =
1$ and $\xi(0) = \xi(1) = \frac{1}{2}$.
\end{thebibliography}
\end{document}